\documentclass[3p,preprint]{elsarticle}
\usepackage{amsthm}
\usepackage{graphicx}
\usepackage{amssymb}
\usepackage{url}
\usepackage{algorithm}
\usepackage{algpseudocode}
\usepackage{amsmath}
\usepackage{bbm}
\DeclareMathOperator*{\avg}{avg}
\newtheorem{definition}{Definition}

\usepackage{color}

\journal{Applied Mathematics and Computation}

\begin{document}

\begin{frontmatter}

\title{Node-Bound Communities for Partition of Unity Interpolation on Graphs}

\author[address-TO,address-GNCS]{Roberto Cavoretto\corref{corrauthor}} 
\ead{roberto.cavoretto@unito.it}

\author[address-TO,address-GNCS]{Alessandra De Rossi}
\ead{alessandra.derossi@unito.it}

\author[address-TO,address-GNCS]{Sandro Lancellotti}
\ead{sandro.lancellotti@unito.it}

\author[address-TO]{Federico Romaniello}
\ead{federico.romaniello@unito.it}

\cortext[corrauthor]{Corresponding author}

\address[address-TO]{Department of Mathematics \lq\lq Giuseppe Peano\rq\rq, University of Torino, via Carlo Alberto 10, 10123 Torino, Italy}

\address[address-GNCS]{Member of the INdAM Research group GNCS}

\begin{abstract}
Graph signal processing benefits significantly from the direct and highly adaptable supplementary techniques offered by partition of unity methods (PUMs) on graphs. In our approach, we demonstrate the generation of a partition of unity solely based on the underlying graph structure, employing an algorithm that relies exclusively on centrality measures and modularity, without requiring the input of the number of subdomains. Subsequently, we integrate PUMs with a local graph basis function (GBF) approximation method to develop cost-effective global interpolation schemes. We also discuss numerical experiments conducted on both synthetic and real datasets to assess the performance of this presented technique.
\end{abstract}

\begin{keyword}
Partition of unity methods \sep Kernel-based approximation \sep graph basis functions \sep graph signal processing \sep graph theory.

\end{keyword}

\end{frontmatter}



\section{Introduction}
Graph signal processing is a very popular tool when it comes to study the field of graph signals and, for this reason, several mathematical techniques such as filtering, compression, noise removal, sampling or decomposition are deeply investigated in the literature, see \cite{ort18,sta19}. Moreover, graphs naturally appear in several problems, as in social networks or traffic maps \cite{cde}. However, they generally have a complex structure that requires efficient and fast processing tools for their analysis. Among the notable characteristics of graphs that depict real systems, the arrangement of vertices into communities stands out. These communities exhibit a high density of edges connecting vertices within them, while the connections between vertices from distinct communities are relatively sparse. Detecting communities is of great importance in sociology, biology and computer science, disciplines where systems are often represented as graphs, nonetheless this problem is very hard and not yet satisfactorily solved, see \cite{for}.

Partition of unity methods (PUMs) enables the execution of tasks like signal reconstruction from samples, node classification, or localized signal filtering on smaller segments of the graph \cite{cde}. Subsequently, they facilitate the reconstruction of the global signal from these localized components.
In the scattered data approximation, the blend of radial basis functions (RBFs) with PUMs results in markedly sparser system matrices in collocation or interpolation tasks, leading to a substantial reduction in computing time, see e.g. \cite{cav21,cav20,fas07}. In what follows, we will focus on approximation methods based on generalized translates of a graph basis function (GBF), see \cite{erb1} and \cite{erb2} for more details.

In this article, we provide a more adaptive technique for the selection of the partitions on the graph than the one presented in \cite{cde}. In particular, by interpreting the communities of a graph as the candidate partitions for the PUM, we propose a process that automatically finds an optimal number of communities or  subdomains based on their modularity and on the underlying graph structure. We remark that the desired number of subdomains are not needed as an input for the GBF-PUM algorithm and, as far as the authors' knowledge, such an approach has never been used before. Our study has been supported by several numerical experiments carried out to validate our new algorithm on both test and real-world datasets.

The structure of the paper is outlined in the following manner. In Section \ref{preli} we give a brief outline of graph theory and positive definite GBFs for signal approximation on graphs. A description of our algorithm, how it can be integrated to the GBF-PUM schema proposed in \cite{cde} is provided in Section \ref{tbsd}, together with a computational complexity analysis. In Section \ref{ne} numerical examples showing the behaviour of the presented technique both in synthetic and real datasets are presented. Finally, Section \ref{concl} concludes the paper also including future work. 

\section{Preliminaries}\label{preli}
\subsection{Graph theory}
A \emph{graph} is an incidence structure $G = (V, E)$, where $V(G)$ is a set whose elements are called \emph{vertices} or \emph{nodes}, and $E(G)$ is a set of pairs of vertices, whose elements are called \emph{edges} or \emph{links}. 

In what follows, we use the letter $G$ to denote a graph and the symbol $|V(G)|$ to denote the cardinality of its vertex set $V(G)$. Moreover, all graphs considered here are connected and simple, that is, without loops, i.e. edges of the form $(v,v)$, and multiple edges, i.e. more than one edge for a given pair of vertices. We refer the reader to \cite{BM} and \cite{newman} for further definitions and notation not explicitly stated in this section.

\begin{definition}
A graph $H=(V(H),E(H))$ is called a \emph{subgraph} of $G=(V(G),E(G))$ if and only if $V(H) \subseteq V(G)$ and $E(H) \subseteq E(G)$. A collection of disjoint subgraphs of $G$ which covers $V(G)$ is a \emph{partition} of $G$ and the subgraphs are called \emph{communities}.    
\end{definition}
Given $u \in V(G)$, the set $N(u)=\left\{v\in V(G)| (u,v) \in E(G)\right\}$ is the \emph{neighbourhood} of $u$. The number of edges $|N(u)|$ incident to a vertex $u$ is said to be the \emph{degree} of $u$, and will be denoted as $deg_G(u)$.    
Given a graph $G$ with $V(G)=\{u_1,\dots,u_n\}$, it is possible to describe its entire structure with matrices, as shown below.
\begin{definition}
 The \emph{adjacency matrix} of the graph $G$ is a square $n \times n$ matrix $A$ such that its entry $A_{ij}=1$ if and only if there is an edge between vertices $u_i$ and $u_j$, and zero otherwise.   
\end{definition}
\begin{definition}
  The \emph{Laplacian matrix} of the graph $G$ is a square $n \times n$ matrix $L$ defined as $L=D-A$, where $D$ is the $n \times n$ diagonal matrix with $D_{ii}=deg_G(u_i)$ and $A$ is the adjacency matrix of $G$.   
\end{definition}
 
Centrality is a very important notion when it comes to identifying the most important nodes in a graph. It may come in different flavours and each one of them highlights its importance from a different perspective. The following definition will be useful.
\begin{definition}
\label{katz}
\emph{Katz centrality} computes the relative influence of a vertex within a graph by measuring the number of the immediate neighbours and also all other vertices in the network that connect to the vertex under consideration through these immediate neighbours. Connections made with distant neighbours are, however, penalized by an attenuation factor $\alpha$, usually set as $0.5$. The Katz centrality for vertex $v_i$ is: 
$$C_{Katz}(u_i)= \sum_{k = 1}^{\infty}  \sum_{j = 1}^{n}  \alpha^k (A^k)_{ij}.$$    
\end{definition}
Furthermore, we will need the following concepts, to better understand the algorithms presented in Section \ref{tbsd}.
\begin{definition}
\label{cut}
Given a graph $G$, a partition of $V(G)$ into two subgraphs $c_1,c_2$ is said to be a \emph{cut} $c=(c_1,c_2)$. The set of edges having one endpoint in $c_1$ and the other one in $c_2$ is said the \emph{cut-set} of $c$. When $s$ and $v$ are specified vertices of $V(G)$, then an $(s,v)$-cut is a cut $c=(c_1,c_2)$ such that $s\in c_1$ and $t \in c_2$.
\end{definition}
\begin{definition}
\label{capacity}
Given an undirected and unweighted graph, the \emph{capacity} of a cut $c$ is the number of edges belonging to the cut-set of $c$. If the graph has weights on the edge-set, the capacity is the sum of the weights of the edges belonging to the cut-set. Accordingly, in an undirected graph, the cut whose weight is smaller (or equal) than the weight of any other cut is called \emph{minimum cut}.
\end{definition}
\begin{definition}
\label{jacard}
Let $u,v \in V(G)$, the \emph{Jaccard similarity index} $J(u,v)= \frac{|N(u)\cap N(v)|}{|N(u)\cup N(v)|}$ measures the similarity between two vertices. It can easily be extended to communities $U,V$ of $G$ as follows: $$J(U,V)=\avg\limits_{u\in U, v \in V}\left\{J(u,v)\right\}.$$ Note that, by design, $0\leq J(U,V) \leq 1$.    
\end{definition}
\begin{definition}
The \emph{modularity} $Q$ of a graph is a measure of the strength of the division of a graph into communities $c_1,\dots, c_n$. Graphs with high modularity have dense connections between the vertices within communities but sparse connections between vertices in different communities. It is defined as 
$$Q=\frac{1}{2m}\sum_{i,j}\left(A_{ij}-\frac{k_ik_j}{2m}\right)\delta(c_i,c_j),$$
where $m$ is the number of edges,  $k_i$ is the degree of $u_i$ and $\delta(c_i,c_j)$ is 1 if $u_i$ and $u_j$ are in the same community, zero otherwise. The value for the modularity for unweighted and undirected graphs lies in the range $[-1/2,1]$ and it is positive if the number of edges within communities exceeds the number expected on the basis of choice.
\end{definition}

\subsection{Signal approximation on graphs with positive definite GBFs}\label{pdgbf}
GBFs, when only a limited number of samples are available, serve as straightforward and efficient instruments for interpolating and approximating graph signals. In particular, when positive definite GBFs are used, better results can be achieved. The theory closely parallels a related one in scattered data approximation using positive definite RBFs in Euclidean space, see \cite{SW} and \cite{Wen}.

\begin{definition}
Given a graph $G$, a function $x: V(G) \rightarrow \mathbb{R}$ is said a \emph{graph signal}. If $|V(G)|=n$, then the $n$-dimensional vector space of all graph signals will be denoted by $\mathcal{L} (G)$. 
\end{definition}


By means of the Laplacian matrix $L$, it is possible to consider a Fourier transform on the graph $G$ since $L=UM_{\lambda}U^T$, where $M_{\lambda}=diag(\lambda_1,\dots,\lambda_n)$ identifies the diagonal matrix whose entries are the (increasingly ordered) eigenvalues of $L$ and the columns of the orthonormal matrix $U$ denote the normalised eigenvectors of $L$ w.r.t. the corresponding eigenvalue. The set $\hat{G}=\{ u_1,\dots,u_n \}$ consisting of the ordered eigenvectors forms an orthonormal basis for the space $\mathcal{L}(G)$ and it is called the \emph{spectrum} of $G$.
\begin{definition}
 Given a graph $G$, let $L=UM_{\lambda}U^T$ be its Laplacian and let $x \in \mathcal{L}(G)$ be a signal on $G$. The Fourier transform on the graph $G$ of the signal $x$ is the vector $\hat{x}=U^Tx$, whose inverse Fourier transform is defined by $x=U\hat{x}$.    
\end{definition}
We will now focus on kernels $K: V\times V \rightarrow \mathbb{R}$ defined on $G$ such that, for every $u,v \in V: K(u,v)=K(v,u)$. Kernels permit to consider a linear operator $\textbf{K}:\mathcal{L}(G) \rightarrow \mathcal{L}(G)$ which acts on $x \in \mathcal{L}(G)$ by $$ \textbf{K}x(v_i)=\sum_{j=1}^n K(v_i,v_j)x(v_j). $$
Hence, $\textbf{K}$ can be represented as the following symmetric  matrix of size $n \times n$:
$$K=\begin{bmatrix}
K(v_1,v_1) & K(v_1,v_2) & \dots & K(v_1,v_n) \\
\vdots & \vdots & \ddots & \vdots \\
K(v_n,v_1) & K(v_n,v_2) & \dots & K(v_n,v_n)
\end{bmatrix}.$$
We will name a symmetric kernel function $K$ \emph{positive definite} if the associated matrix $K$ is strictly positive definite, that is $x^TKx \geq 0$, $\forall x \in \mathbb{R}^n, x \neq 0$. We remark that every positive definite kernel $K$ induces on the space $\mathcal{L}(G)$ an inner product $\langle x,y \rangle_{K}=y^TKx$ and a norm $\|x\|_{K}$ , for every $x,y \in \mathcal{L}(G)$. This resulting inner product space will be referred to as \emph{native space} $\mathcal{N}_{K}$. Moreover, it is a reproducing kernel Hilbert space, see \cite{Aro}. To obtain a GBF approximation $x_{\star}$ for a signal $x$ based on $N$ samples $x(w_i)$, where $i=1,\dots,N$, we seek the solution $x_{\star}$ to the following regularized least squares (RLS) problem:
\begin{equation}
  x_{\star}=\text{argmin}_{y \in \mathcal{N}_{\textbf{K}}}\left\{\frac{1}{N}\sum_{i=1}^N|x(w_i)-y(w_i)|^2+\gamma \|y\|^2_{\textbf{K}_f} \right\}.  \label{RLSprob}
\end{equation}
The first term in of Equation (\ref{RLSprob}) is referred to as the data fidelity term, which guarantees the proximity of the values $x_{\star}(w_i)$ to the provided values $x(w_i)$, on the
sampling set $W$, whereas $\gamma >0$ is denoted as the regularization parameter.

 Numerical examples presented in Section \ref{ne} rely on polyharmonic splines on graphs, which are an example of a positive definite GBF based on the kernel:
$$K_{f_{(\epsilon I_n+L)^{-s}}}=(\epsilon I_n+L)^{-s}=\sum_{k=1}^{n} \frac{1}{(\epsilon + \lambda_k)^s}u_ku_k^T,$$
where $u_1,\dots,u_n$ are the eigenvectors of the Laplacian $L$. If $\lambda_1$ is the largest eigenvalue of $L$, this kernel is positive definite for $\epsilon > -\lambda_1$ and $s > 0$ (see \cite{pes} and \cite{wanawa}). Hence values of $\epsilon= -\lambda_1 +1$ and $s=1$ are set for the experiments.

\section{Topology based subdomain detection}\label{tbsd}
Signal interpolation on graphs with PUM requires overlapping subdomains with the property that at least one vertex of the subdomain is an interpolation node. In what follows we describe how it is possible to find such subdomains for PUM on graphs, based on the structural topology of the underlying structure.
\subsection{Overlapping communities detection}
Community detection is an edge-cutting topic in graph theory, indeed different approaches may be used, see \cite{newman} for a detailed overview. Our method is a divisive technique where the centrality of the interpolation nodes is the main criterion for splitting, allowing overlapping among communities. 
The finding communities process was decomposed into five different functions combined with each other. It is worth noting here that, when the adjacency matrix $A$ of the graph $G$ is part of the input of an algorithm, automatically all information on nodes $V(G)$ and edges $E(G)$ is stored. In detail, the main function is Algorithm \ref{community_detection}, in which, given in input $A$, the adjacency matrix of the graph $G$, and the set of the interpolation (or sample) nodes $W$, an attempt is made to iteratively divide communities until the number of interpolation nodes is reached or no attempt at division leads to an increase in modularity. 
During the attempt of division, Algorithm \ref{find_partition}  provides a new candidate sub-partition, for each community, with its modularity obtained by splitting the considered community. In the end, the sub-partition with the highest modularity, if greater than the initial modularity, superseded the previous partition.
Then, small communities are merged with the big ones using Algorithm \ref{join_communities} and, finally, with Algorithm \ref{expand_communities}, communities are expanded to create the required overlapping for the GBF-PUM.

\begin{algorithm}
    \caption{community\_detection}
    \hspace*{\algorithmicindent} \textbf{Input:}
    \vspace{-3mm}
    \begin{itemize}
        \setlength\itemsep{-0.2em}
        \item[] $A$: adjacency matrix of the graph $G$, $W$: interpolation vertices, $r$: ratio of neighbours of a vertex inside the community over the total number, $dmax$: distance of max augmentation, $dmin$: distance of min augmentation.  
    \end{itemize}
    
    \begin{algorithmic}
        \State $partition \rightarrow \lbrace V(G) \rbrace$
        \State $Q' \rightarrow -1$
        \State $Q\rightarrow modularity(partition)$
        \While{$|partition| <= |W| \quad \& \quad  Q' < Q $}
            \State $Q'\rightarrow Q$
            \For{$i \leq |partition|$}    
                \State $partition_i, Q_i \rightarrow find\_partition(A, partition, partition(i), W)$
                \If{$Q_i > Q$}
                    \State $partition \rightarrow partition_i$
                    \State $Q \rightarrow Q_i$
                \EndIf
            \EndFor
        \EndWhile
        \State $partition \rightarrow join\_communities(A, partition)$
        \State $expanded\_partition \rightarrow expand\_communities(A, partition, r, dmax, dmin)$
    \end{algorithmic}
    \hspace*{\algorithmicindent} \textbf{Output:}
        \vspace{-3mm}
\begin{itemize}
\setlength\itemsep{-0.2em}
    \item[] $partition$: partition of $V$,
    $expanded\_partition$: expanded partition with overlap of $V$.
    \end{itemize}
\label{community_detection}
\end{algorithm}

Algorithm \ref{find_partition} is used for generating the candidate sub-partition. It takes in input the adjacency matrix $A$ of the graph $G$, a $partition$ of $V(G)$, the community to divide and the set of samples $W$. It begins calling  Algorithm \ref{split_net} to retrieve the two sub-communities obtained by the splitting, then it generates the sub-partition and evaluates its modularity, returning them as output.

\begin{algorithm}
    \caption{find\_partition}
    \hspace*{\algorithmicindent} \textbf{Input:}
    \vspace{-3mm}
    \begin{itemize}
        \setlength\itemsep{-0.2em}
        \item[] $A$: adjacency matrix of the graph $G$, $partition$: partition of $V(G)$, $community$: a community of $G$,  $W$: interpolation vertices.   
    \end{itemize}
    
    \begin{algorithmic}
        \State $c_1$, $c_2$ $\rightarrow split\_net(A, community, W)$
        \State $partition' \rightarrow (partition \setminus community) \cup \lbrace c_1, c_2 \rbrace$
        \State $Q' \rightarrow$ modularity of $partition'$
    \end{algorithmic}
    \hspace*{\algorithmicindent} \textbf{Output:}
        \vspace{-3mm}
\begin{itemize}
\setlength\itemsep{-0.2em}
    \item[]$partition'$: the sub-partition, $Q'$: modularity of $partition'$.
    \end{itemize}
\label{find_partition}
\end{algorithm}

The core step of the entire procedure is mastered by Algorithm \ref{split_net}, which, taken as input a connected graph (which may be $G$ or one of its communities) and the set of interpolation nodes $W$, evaluates the Katz centrality for each sample vertex (see Definition \ref{katz}). After that, it selects the two vertices $s$ and $v$ with the highest centrality values. Furthermore, the capacity (see Definition \ref{capacity}) of all edges having as first end-point $s$ or $v$ and as the other end-point neither a vertex different from $v$ or $s$ nor one of their neighbours respectively is set to infinity; whereas, a capacity value of 1 is set to all other edges.
Eventually, it performs a minimum $(s, v)$-cut (see Definition \ref{cut}) obtaining the 2 sub-communities returned as output.

\begin{algorithm}
    \caption{split\_net}
    \hspace*{\algorithmicindent} \textbf{Input:}
    \vspace{-3mm}
    \begin{itemize}
        \setlength\itemsep{-0.2em}
        \item[] $A$: adjacency matrix of the graph $G$, $G$: connected graph $G$, $W$: interpolation vertices.   
    \end{itemize}
    
    \begin{algorithmic}
        \State Evaluate the Katz centrality for each vertex in $W$
        \State Select $s$ and $v$, the vertices with high centrality
        \For{$\{x, y\} \in E(G)$}
            \State $capacity(\{x, y\}) \rightarrow 1$
        \EndFor
        \For{$u \in N(s)$}
            \If{ $u \notin N(v) \cup \lbrace v \rbrace$}
                \State $capacity(\lbrace s, u \rbrace) \rightarrow \infty$
            \EndIf
        \EndFor
        \For{$u \in N(v)$}
            \If{ $u \notin N(s) \cup \lbrace s \rbrace$}
                \State $capacity(\lbrace v, u \rbrace) \rightarrow \infty$
            \EndIf
        \EndFor
        \State $c_1, c_2 \rightarrow \text{Compute the value and the node partition of a minimum $(s, v)-cut$}$
    \end{algorithmic}
    \hspace*{\algorithmicindent} \textbf{Output:}
        \vspace{-3mm}
\begin{itemize}
\setlength\itemsep{-0.2em}
    \item[] $c_1$: split community,  $c_2$: split community, $partition'$: new partition, $Q'$: modularity of $partition'$.
    \end{itemize}
\label{split_net}
\end{algorithm}

Algorithm \ref{join_communities} merges the small communities into the most similar large one. Taking as input $A$ and the $partition$ it starts by dividing the communities into small if contains less than the $2\%$ of vertices and large. For each small community, it evaluates the Jaccard similarity (see Definition \ref{jacard}) among the small community considered and all the large ones in the graph. Then, the small one is merged with the most similar big community.

\begin{algorithm}
    \caption{join\_communities}
    \hspace*{\algorithmicindent} \textbf{Input:}
    \vspace{-3mm}
    \begin{itemize}
        \setlength\itemsep{-0.2em}
        \item[] $A$: adjacency matrix of the graph $G$, $partition$: partition of $V(G)$.   
    \end{itemize}
    
    \begin{algorithmic}
        \State Divide $partition$ in $small\_communities$ and $big\_communities$ (a community is small if contains less than the $2\%$ of vertices) 
        \For{$community$ in $small\_communities$}
            \State find the highest Jacard index between the $community$ and the $big\_communities$
            \State merge the communities with the highest Jacard index
            \State update $partition$
        \EndFor
    \end{algorithmic}
    \hspace*{\algorithmicindent} \textbf{Output:}
        \vspace{-3mm}
\begin{itemize}
\setlength\itemsep{-0.2em}
    \item[] $partition$: partition of $V$ without small communities.
    \end{itemize}
\label{join_communities}
\end{algorithm}

The last step in the search for communities is completed by Algorithm \ref{expand_communities}, which takes in input $A$, $partition$, $r$, $dmax$ and $dmin$. Then, for all vertices in a community it counts the ratio of the number of neighbours inside the community among the total number of neighbours. If this ratio is less than $r$, it expands the community with all vertices within 
 distance $dmax$ from the considered vertex, otherwise the community is expanded with all vertices within distance $dmin$. The overlapping communities needed for the GBF-PUM interpolation (Algorithm \ref{gbf_pum}) are given as output.
\begin{algorithm}
    \caption{expand\_communities}
    \hspace*{\algorithmicindent} \textbf{Input:}
    \vspace{-3mm}
    \begin{itemize}
        \setlength\itemsep{-0.2em}
        \item[] $A$: adjacency matrix of the graph $G$, $partition$: partition of $V(G)$, $r$: ratio of neighbours of a vertex inside the community over the total number, $dmax$: distance of max augmentation, $dmin$: distance of min augmentation.   
    \end{itemize}

    \begin{algorithmic}
        \For{$i \leq |partition|$}
            \For{$v$ in $partition(i)$}
                \State $N(v) \rightarrow$ neighbours of $v$
                \State $N(v)_{in} \rightarrow$ neighbours of $v$ inside $community$
                \If{$|N(v)_{in}| < r*|N(v)|$}
                    \State add to $partition(i)$ all the vertices within distance of $dmax$ from $v$
                \Else
                    \State add to $partition(i)$ all the vertices within distance of $dmin$ from $v$
                \EndIf
            \EndFor
        \EndFor
    \end{algorithmic}
    \hspace*{\algorithmicindent} \textbf{Output:}
        \vspace{-3mm}
\begin{itemize}
\setlength\itemsep{-0.2em}
    \item[] $partition$: partition of $V$ with expanded communities.
    \end{itemize}
\label{expand_communities}
\end{algorithm}


\subsection{Computational complexity of Algorithm \ref{community_detection}}
To analyse the computational complexity of Algorithm \ref{community_detection} we have to evaluate the complexity of the other algorithms called during the procedure. We remark here that logarithms are taken base 2.
In details, in Algorithm \ref{split_net}, \emph{split\_net}, the minimum-cut procedure costs $O(|c_i||E(c_i)|+|c_i|^2\log(|c_i|)$ for each community $c_i$. Hence its total computation time in one iteration is $$\sum_{i=1}^{|C|}O(|c_i||E(c_i)|+|c_i|^2\log(|c_i|) \leq O(|V||E|+|V|^2\log(|V|).$$
On the other hand, calculating the Katz centrality for the vertices in $G$ takes $O(|V|+|E|)$, when the power iteration algorithm is used. Hence, Algorithm \ref{split_net} has the same complexity as the minimum-cut.  

The computation of the modularity for each community should also be considered when estimating the total complexity. As the changes of $Q$ are tracked during the execution, its complexity is bounded by $O(|V|+|E|)$. Moreover, calculating the Jaccard index between the communities takes $O(|V|^2)$ in the worst-case scenario $|C|=|V|$. All the other executions of the procedures take linear time. Hence, both the evaluation of modularity and Jaccard index do not affect the total time complexity as they are faster than Algorithm \ref{split_net}. We can conclude that in the worst possible case when the number of iterations is exactly equal to $|V|$, the total time complexity is $$O(|V|^2|E|+|V|^3\log(|V|)\approx O(|V|^3\log(|V|),$$
which is in polynomial time.
\subsection{GBF-PUM approximation on graphs}
Partition of unity is a widespread tool used in mesh-free approximation problems as it sharply reduces the computational cost when related subproblems are solved \cite{fas15}. In \cite{cav22}  and \cite{cde} this technique has been adapted to graph signal interpolation and it may be applied to social networks or traffic map problems because of their underlying graph structure.
Once the overlapping communities of a graph are found, to construct a partition of unity we need a set of functions $\varphi^{(j)} \in \mathcal{L}(G)$, $j=1,\dots, |C|$, such that:
\begin{itemize}
    \item $\rm{supp}(\varphi^{(j)})\subseteq V_j$;
    \item $\varphi^{(j)} \geq 0$;
    \item $\sum_{j=1}^{|C|}\varphi^{(j)}(v)=1$, $\forall v \in V$.
\end{itemize}
In particular, we choose the characteristic function $\varphi^{(j)}(v)=\mathbbm{1}_{V_{j}}(v)=\begin{cases}
1\text{ if $v \in V_j$}\\
0\text{ if $v \not\in V_j$}
\end{cases}$,  to get the following partition of unity:
$$ \varphi^{(j)}(v)=\frac{\mathbbm{1}_{V_{j}}(v)}{|j \in \{1,\dots,|C|:v  \in V_j\}|}.$$
The main step of this process is the calculation of the local approximants $x_{\star}^{(j)}$ on the subgraph $G_j$ corresponding to the $j$-th community. This is done by using a GBF approximation scheme, as shown in Subsection \ref{pdgbf}. When the local approximants are computed, the global GBF-PUM approximation scheme on the entire graph $G$ is:
$$x_{\star}(v)=\sum_{j=1}^{|C|}\varphi^{(j)}(v)x_{\star}^{(j)}(v).$$
We remark that overlapping subdomains are needed for a good PUM. 
Further details on the error of the global approximant can be found in \cite{cde}.
The following Algorithm \ref{gbf_pum} traces the steps to obtain the global interpolant. Unlike the other pseudocodes presented above, Algorithm \ref{gbf_pum} is described verbosely as all other procedures except community detection are well-known in literature. The reader may find further details about them in \cite[Section 4]{cde}.



\begin{algorithm}
    \caption{GBF-PUM Approximation on Graphs}
    \hspace*{\algorithmicindent} \textbf{Input:}
    \vspace{-3mm}
    \begin{itemize}
        \setlength\itemsep{-0.2em}
        \item[] $A$: adjacency matrix of the graph $G$, $W$: interpolation vertices, $r$: ratio of neighbours of a vertex inside the community over the total number, $dmax$: distance of max augmentation, $dmin$: distance of min augmentation.  
    \end{itemize}
    
    \begin{algorithmic}
    
\State Find the overlapping communities $C$ applying community detection (Algorithm \ref{community_detection})
\State Construct a partition of unity subordinate to the communities in $C$
\State For all subgraphs $c$ in $C$ calculate the local Laplacian $L^c$
\State For all subgraphs $c$ in $C$ construct the local GBF kernel
\State For all subgraphs $c$ in $C$ calculate the local GBF approximant
\State Create a global GBF-PUM approximation from the local ones    
\end{algorithmic}
    \hspace*{\algorithmicindent} \textbf{Output:}
        \vspace{-3mm}
\begin{itemize}
\setlength\itemsep{-0.2em}
    \item[] $x_{\star}$: GBF-PUM approximant of the signal.
    \end{itemize}
\label{gbf_pum}
\end{algorithm}


\section{Numerical examples}\label{ne}

In what follows, we exhibit the subdivision ability of the algorithm by considering Zachary's Karate Club graph \cite{Zak}. Then, numerical results obtained by applying our framework on the Minnesota Road Graph \cite{minnroadgraph}, using a test function as a signal on its vertices, and on the real dataset of the freight transport dataset of Brussels \cite{brunet} are presented. To perform tests, the values of $r=0.75$, $dmax=6$, $dmin=4$ are set as input for Algorithm 1.

All tests were carried out on the infrastructure for high-performance computing \emph{MathHPC}, virtual cloud server of the structure \emph{HPC4AI} (High-Performance Computing for Artificial Intelligence: \url{https://www.dipmatematica.unito.it/do/progetti.pl/View?doc=Laboratori_di_ricerca.html}).

\subsection{Zachary's Karate Club graph}
\begin{figure}[h!]
    \centering
    \includegraphics[width=.5\textwidth]{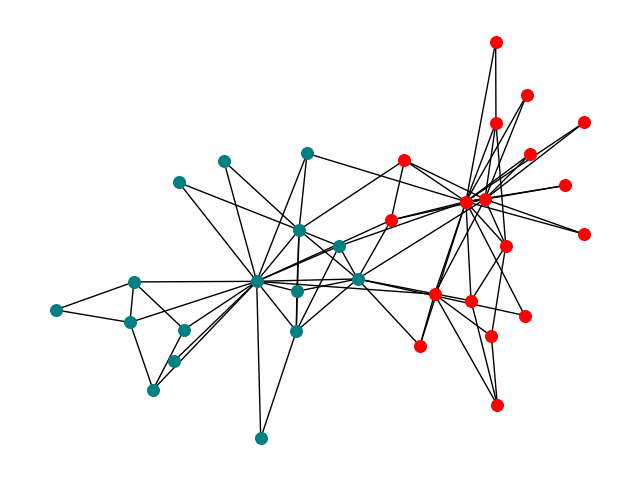} 
    \caption{The Zachary's Karate Club graph two communities.}
    \label{kt_comm}
\end{figure}
\begin{figure}[h!]
    \centering
    \includegraphics[width=.4\textwidth]{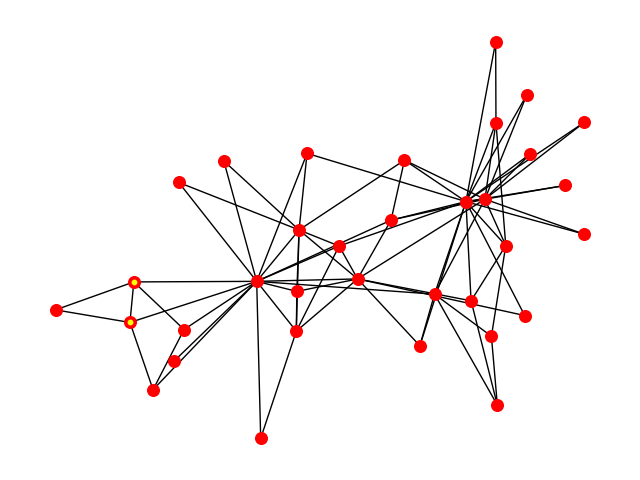} 
\includegraphics[width=.4\textwidth]{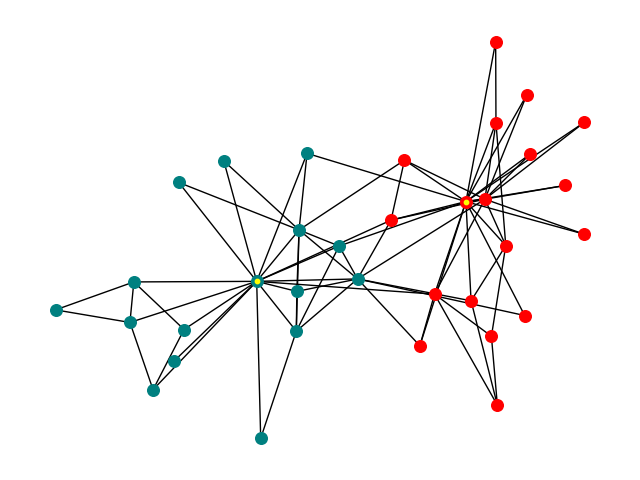}
    \caption{On the left it is shown how a bad selection of the sample nodes leads to a non-subdivision of the network. On the right, when the two sample nodes are chosen as the instructor and the administrator of the karate club the algorithm finds the ground truth two communities. The sample nodes are dotted in yellow. }
    \label{badgoodcomm}
\end{figure}
Zachary's Karate Club graph is a 34-vertex, 78-edge social network of a University Karate Club, described in the paper \cite{Zak}. During the study, a conflict arose between the administrator and the instructor, hence splitting the club into two parts. Half of the members formed a new club around the instructor; members from the other part found a new instructor or gave up karate. These ground-truth communities are depicted in Figure \ref{kt_comm}. 

Let $W=\{u_1,u_2\}$ be the set of sampling nodes for the community algorithm, and suppose that $u_1,u_2$ are neighbours. An execution of Algorithm \ref{community_detection} does not perform any subdivision as there is no way to increase the modularity by splitting the two sample nodes. On the other hand, if we choose $u_1$ and $u_2$ as the node corresponding to the instructor and the administrator, our algorithm performs exactly one subdivision splitting the network into the well-known two communities, see Figure \ref{badgoodcomm}.

\subsection{Minnesota Road Graph}
The Minnesota Road Graph is a 2642-vertex, 3304-edge network representing, respectively, intersections and roads of the State of Minnesota \cite{shu16}. A test function based on the Laplacian matrix introduced in \cite{cav22} and \cite{cde} has been adapted and used as a signal on the vertex set to apply the presented algorithms.
\begin{figure}[h!]
    \centering
    \includegraphics[width=.5\textwidth]{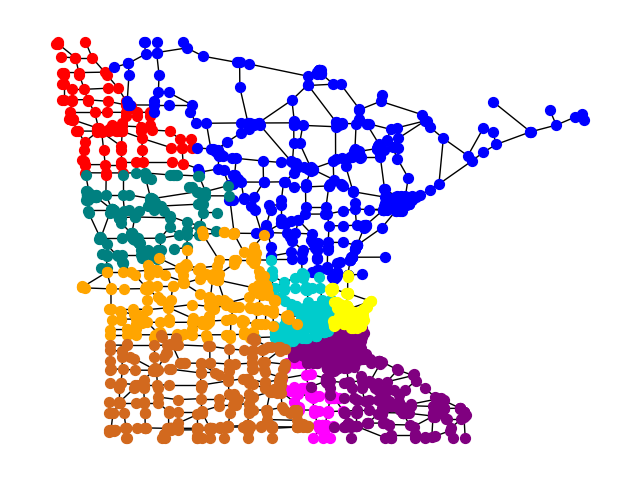} \\
\includegraphics[width=.5\textwidth]{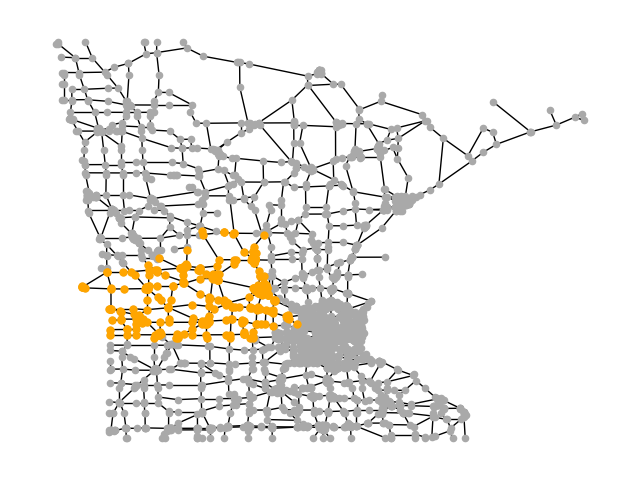}\hfill
\includegraphics[width=.5\textwidth]{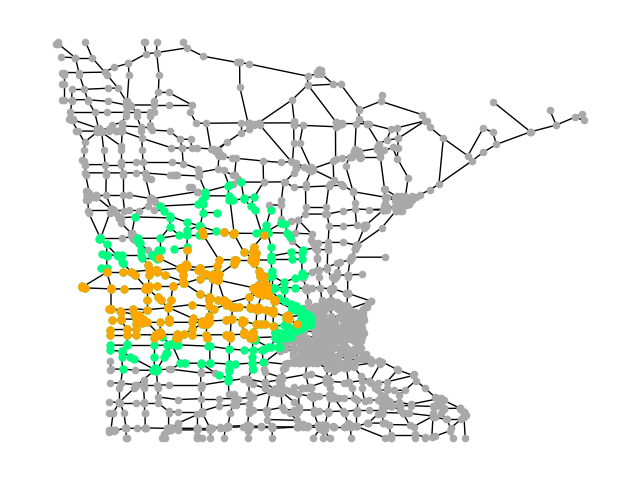}
    \caption{On top, the partition of the vertices in 9 different disjoint subdomains is shown. It is obtained using 800 interpolation nodes. On bottom, a subdomain is highlighted to show how it is expanded when creating the overlap (left to right).}
    \label{subdomain_extension}
\end{figure}
The communities are shown in Figure \ref{subdomain_extension}, whereas numerical results are summarised in Table \ref{tabella_minessota}.
We would like to remark that the number of communities found by the divisive Algorithm \ref{community_detection} only depends on the centrality values of the sample vertices.
The number of detected communities, Relative-Max-Absolute-Error (RMAE),  Relative-Root-Mean-Squared-Error (RRMSE), and computational time to the variation of the number of interpolation vertices are shown in Table \ref{tabella_minessota}.


\begin{table}[h]
    \centering
    \begin{tabular}{|c|c|c|c|c|}
    \hline
    $W$       & Communities & RMAE &    RRMSE &     time (s) \\
    \hline
    400	 & 11	&  3.509e-01 & 3.208681e-02	&  1.425e+02 \\
    \hline
    800	 & 9 &	 9.276e-02 & 5.880306e-03	& 1.879e+02 \\
    \hline
    1200	& 11 & 1.977e-02 &	1.444769e-03	& 1.587e+02 \\
    \hline
    1600	& 10 &	 1.517e-02 & 6.787503e-04	& 1.586e+02  \\
    \hline
    2000     &   9 &	 5.078e-03 & 	2.122356e-04	& 1.690e+02 \\
    \hline
    \end{tabular}
    \caption{A summarisation of the numerical results for increasing number of sampling nodes in the Minnesota Road Graph. The number of communities,  RMAE, RRMSE, and computational time are also shown.}
    \label{tabella_minessota}
\end{table}

\subsection{Brussel Freight Transport network}
\begin{figure}[h]
    \centering
    \includegraphics[width=.4\textwidth]{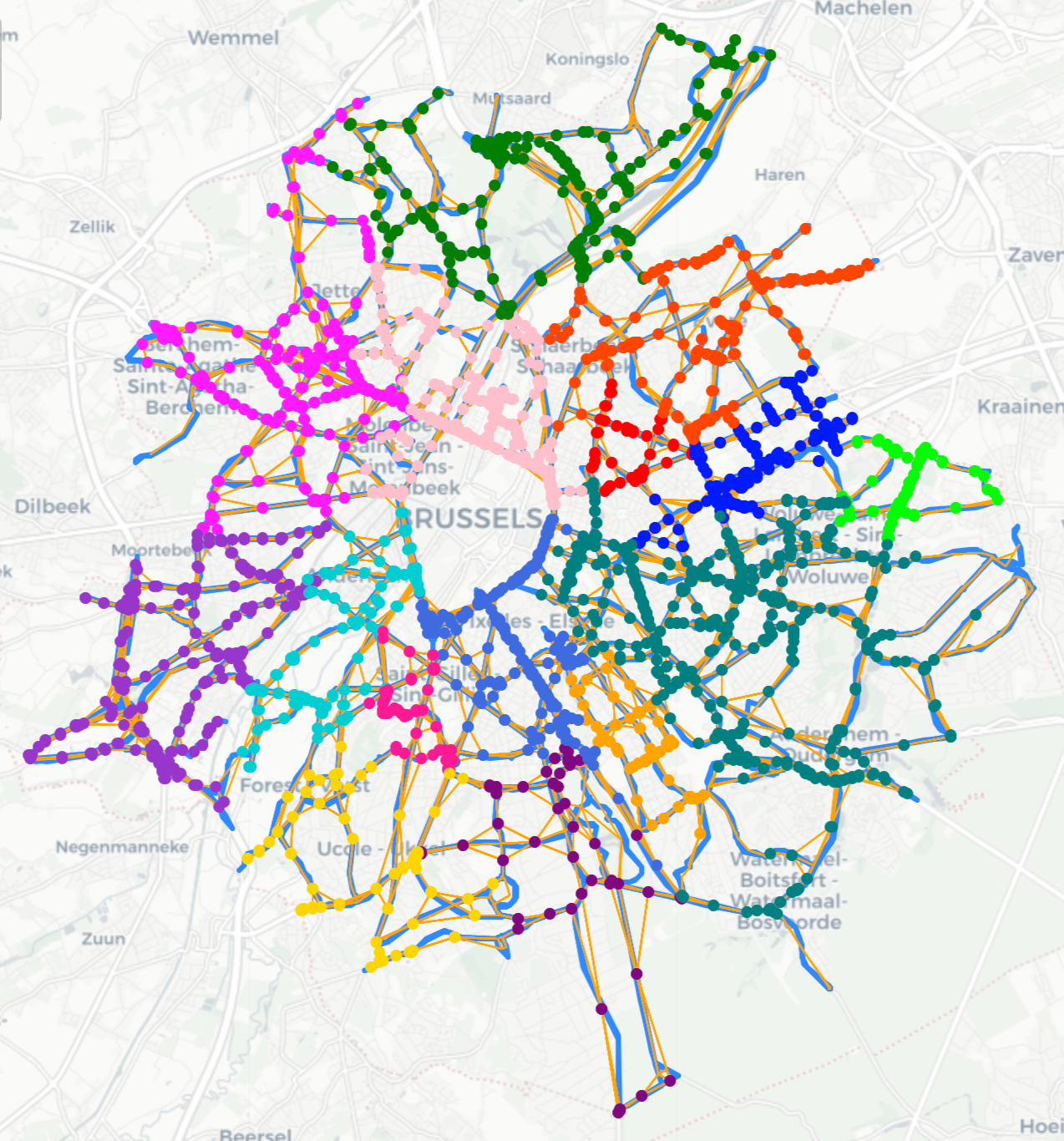} 
\includegraphics[width=.4\textwidth]{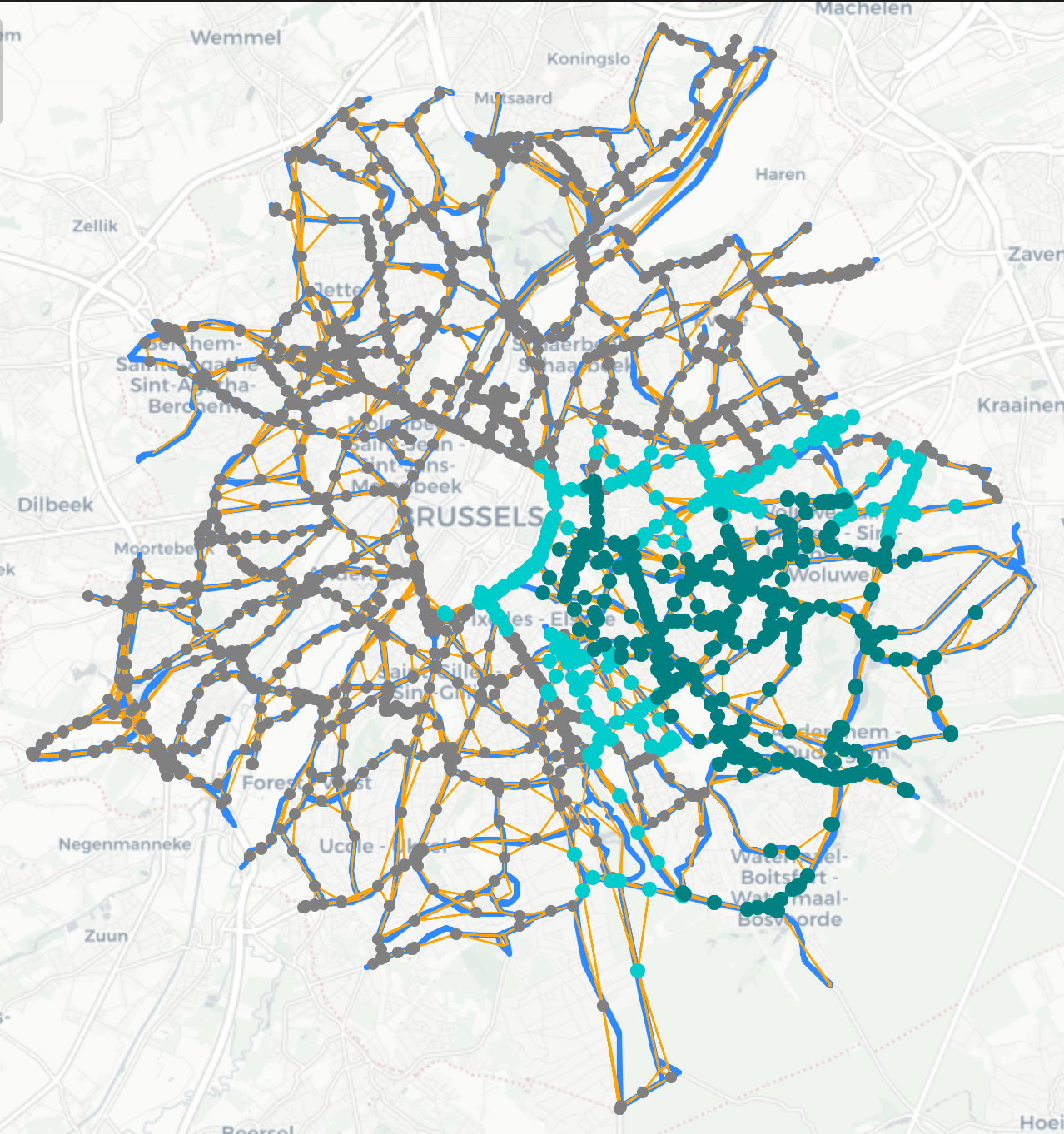} \\
\includegraphics[width=.4\textwidth]{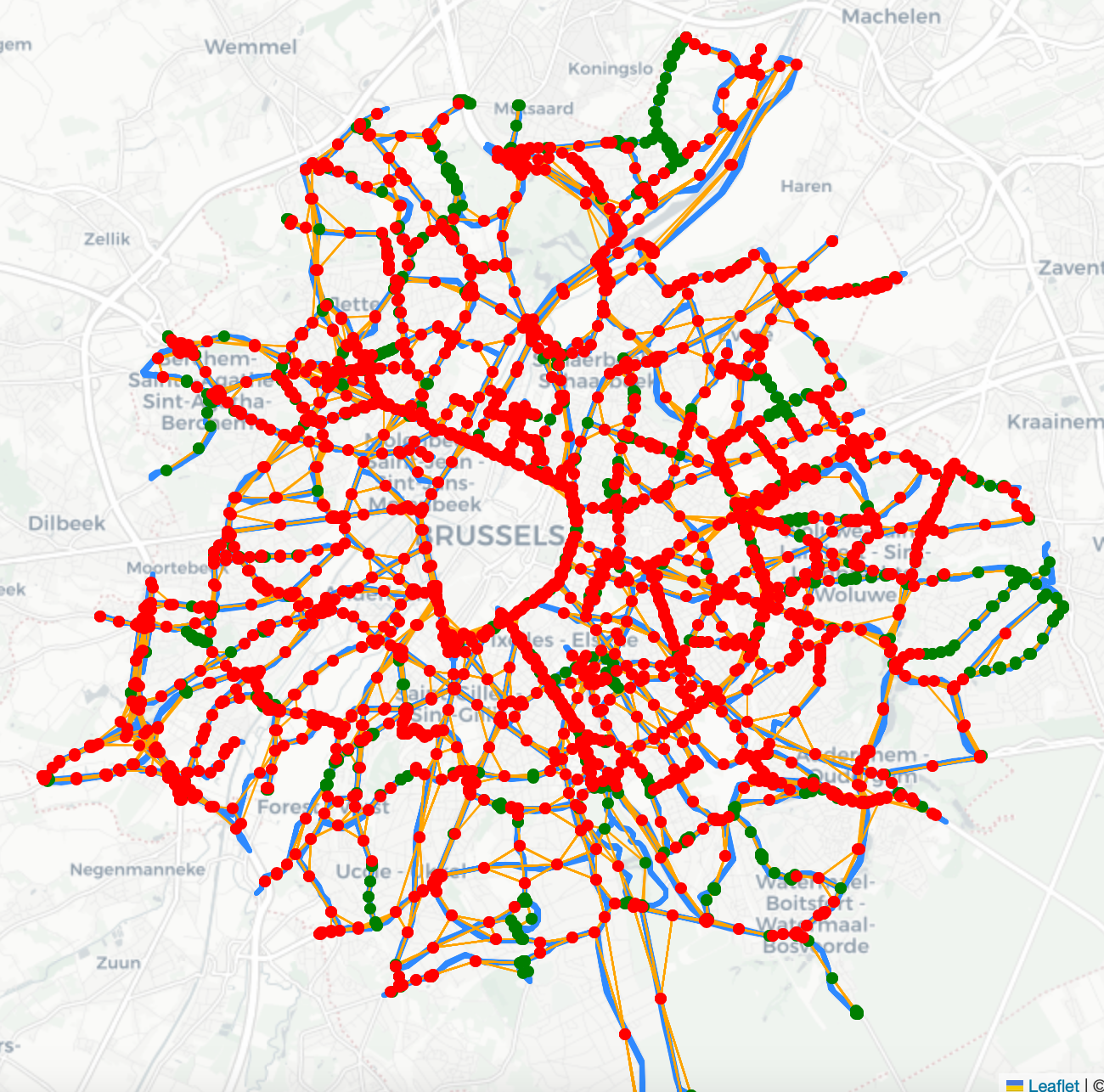}

    \caption{On top the division in communities of the Brussel Freight Transport Network using 1200 sample nodes and the expansion of one of them is shown (left to right). On bottom the sample set of measures is highlighted in red, whereas the unknown values to predict are in green.}
    \label{brux}
\end{figure}

The Brussel Freight Transport Network is a 4540-vertex, 14946-edge graph representing the road network of the city of Brussels \cite{brunet}. For 3451 vertices, a measurement of the flow every 15 minutes is available. As signal for the numerical test, we will consider a single-time measure of the flow in the sample vertices. The associated real data problem will be reconstructing the missing value via GBF-PUM interpolation, whereas for testing purposes we considered the largest connected component of the vertices on which we know the flow. A sample of found communities is shown in Figure \ref{brux}, while numerical results are reported in Table  \ref{tabella_brussels}.

\begin{table}[h]
    \centering
    \begin{tabular}{|c|c|c|c|c|}
    \hline
    $W$       & Communities & RMAE &    RRMSE &     time (s) \\
    \hline
    400	 & 15	&  1.500e+01 & 6.881386e-02	&  2.767e+02 \\
    \hline
    800	 & 7 &	 9.173e-01 & 6.087108e-02	& 2.752e+02 \\
    \hline
    1200	& 15 & 8.015e-01 &	4.673573e-02	& 2.768e+02 \\
    \hline
    1600	& 9 &	 6.286e-01 & 3.925505e-02	& 2.253e+02  \\
    \hline
    2000     &   17 &	 5.845e-01 & 	3.011343e-02	& 2.168e+02 \\
    \hline
    2400	 & 16	& 6.538e-01	 	& 2.435230e-02	 & 2.166e+02\\
    \hline 
    2800	 & 12	& 3.589e-01	 &	1.770847e-02	& 1.974e+02 \\
    \hline 
    3200	& 6 &	 1.732e-01	 &	7.304001e-03	& 1.629e+02	\\
    \hline
    \end{tabular}
    \caption{A summarisation of the numerical results for increasing number of sampling nodes in the Brussel Freight Transport Network. The number of communities, RMAE, RRMSE and computational time are also shown.}
    \label{tabella_brussels}
\end{table}


\section{Conclusions} \label{concl}
From the numerical experiments (Table \ref{tabella_minessota} and Table \ref{tabella_brussels}) we note that when the number of sample nodes increases, also does the precision of the GBF-PUM method. We remark that in the case of a real dataset, this behaviour is less noticeable but still significant, due to the nature of the data. Moreover, the proposed method does not require the number of subdomains for the PUM as an input and it goes in the direction of finding a reasonable number of subdomains for PUMs in an automatic way. This is a well-studied problem not only when working with graphs or networks, but also in the context of interpolation and approximation of sparse data, when solving PDEs \cite{lar17} and within a stochastic framework \cite{cav23}. As future work, an automatic routine also for community augmentation may be investigated.

\subsection*{Acknowledgments}  

This work has been supported by the INdAM--GNCS 2022 Project \lq\lq Computational methods for kernel-based approximation and its applications\rq\rq, code CUP$\_$E55F22000270001, and by the Spoke ``FutureHPC \& BigData'' of the ICSC – National Research Center in "High-Performance Computing, Big Data and Quantum Computing", funded by European Union – NextGenerationEU. Moreover, the work has been supported by the Fondazione CRT, project 2022 \lq\lq Modelli matematici e algoritmi predittivi di intelligenza artificiale per la mobilit$\grave{\text{a}}$ sostenibile\rq\rq.


\begin{thebibliography}{99}



\bibitem{Aro} N. Aronszajn, Theory of reproducing kernels, Trans. Amer. Math. Soc. 68 (1950), 337--404.

\bibitem{BM} J.A. Bondy, U.S.R. Murty, Graph Theory, Springer Series: Graduate Texts in Mathematics. Springer London, 2008.

\bibitem{cav21} R. Cavoretto, Adaptive radial basis function partition of unity interpolation: A bivariate algorithm for unstructured data, J. Sci. Comput. 87 (2021), 41.

\bibitem{cav20} R. Cavoretto, A. De Rossi, Error indicators and refinement strategies for solving Poisson problems through a RBF partition of unity collocation scheme, Appl. Math. Comput. 369 (2020), 124824.

\bibitem{cav22} R. Cavoretto, A. De Rossi, W. Erb, GBFPUM - A MATLAB Package for Partition of Unity Based Signal Interpolation and Approximation on Graphs, Dolomites Res. Notes Approx. 15 (2022), 25--34.

\bibitem{cde} R. Cavoretto, A. De Rossi, W. Erb, Partition of Unity Methods for Signal Processing on Graphs, J. Fourier Anal. Appl. 27 (2021), 66. 

\bibitem{cav23} R. Cavoretto, A. De Rossi, E. Perracchione, Learning with partition of unity-based kriging estimators, Appl. Math. Comput. 448 (2023), 127938.

\bibitem{erb1} W. Erb, Graph Signal Interpolation with Positive Definite Graph Basis Functions, Appl. Comput. Harmon. Anal. 60 (2022), 368--395.

\bibitem{erb2} W. Erb, Semi-Supervised Learning on Graphs with Feature-Augmented Graph Basis Functions. arXiv:2003.07646 (2020).

\bibitem{fas07} G.E. Fasshauer, Meshfree Approximation Methods with MATLAB, World Scientific, Singapore, 2007, 500 pp.

\bibitem{fas15} G.E. Fasshauer, M.J. McCourt, Kernel-based Approximation Methods Using MATLAB, World Scientific, Singapore, 2015, 536 pp.

\bibitem{for} S. Fortunato, Community detection in graphs, Phys. Report 486 (2010), 75--174.

\bibitem{brunet} T. Guns, S. Hadavi, C. Macharis, S. Verlinde, W. Verbekem Monitoring Urban-Freight Transport Based on GPS Trajectories of Heavy-Goods Vehicles, IEEE Transactions on Intelligent Transportation Systems 20(10) (2019), 3747--3758.

\bibitem{lar17} E. Larsson, V. Shcherbakov, A. Heryudono, A least squares radial basis function partition of unity method for solving PDEs, SIAM J. Sci. Comput. 39  (2017), A2538--A2563. 

\bibitem{newman} M. Newman, Networks: An Introduction, OUP Oxford, 2010, 784 pp.

\bibitem{ort18} A. Ortega, P. Frossard, J. Kovačević, J.M.F. Moura, P. Vandergheynst, Graph signal processing: overview, challenges, and applications, Proc. IEEE 106(5) (2018), 808--828.

\bibitem{pes} I.Z. Pesenson, Variational splines and Paley-Wiener spaces on combinatorial graphs, Constr. Approx. 29(1) (2009), 1--21.

\bibitem{minnroadgraph} R.A. Rossi, N.K. Ahmed, The Network Data Repository with Interactive Graph Analytics and Visualization, in Proceedings of the Twenty-Ninth AAAI Conference on Artificial Intelligence, 2015, \emph{http://networkrepository.com}.

\bibitem{SW} R. Schaback, H. Wendland, Approximation by Positive Definite Kernels, in Advanced Problems in Constructive Approximation, Birkhäuser Verlag, Basel, 2003, pp. 203--222


\bibitem{shu16} D.I. Shuman, B. Ricaud, P. Vandergheynst, Vertex-frequency analysis on graphs, Appl. Comput. Harm. Anal. 40(2) (2016), 260--291.

\bibitem{sta19} L. Stankovi{\'c}, L. Dakovi{\'c}, E. Sejdi{\'c}, Introduction to Graph Signal Processing,  in Vertex-Frequency Analysis of Graph Signals, Springer, 2019, pp. 3--108.

\bibitem{wanawa} J.P. Ward, F.J. Narcowich, J.L. Ward, Interpolating splines on graphs for data science applications, Appl. Comput. Harm. Anal. 49(2) (2020), 540--557.

\bibitem{Wen} H. Wendland, Scattered Data Approximation, Cambridge University Press, Cambridge, 2005, 336 pp.

\bibitem{Zak} W.W. Zachary, An Information Flow Model for Conflict and Fission in Small Groups, Jour. of Anthropological Research 33(4) (1977), 452--473.









 



\end{thebibliography}
\end{document}